\input amstex
\documentstyle{amsppt}
\topmatter
\magnification=\magstep1
\pagewidth{5.2 in}
\pageheight{6.7 in}
\abovedisplayskip=12pt \belowdisplayskip=12pt
\NoBlackBoxes
\title A note on $p$-adic invariant integral in the rings of $p$-adic integers\endtitle

\author Taekyun Kim\endauthor
\affil{ {\it Jangjeon Research Institute for Mathematical Sciences and Physics,\\
 252-5 Hapcheon-Dong Hapcheon-Gun Kyungnam, 678-802, S. Korea\\
 e-mail: tkim$\@$kongju.ac.kr/ tkim64$\@$hanmail.net}}
\endaffil

\abstract{In [2], I constructed the $p$-adic $q$-integral $I_q(f)$
on $\Bbb Z_p$. In this paper, we consider the properties of
$\lim_{q\rightarrow -1}I_q(f)=I_{-1}(f).$ Finally we give the some
applications of $I_{-1}(f)$ and integral equations for
$I_{-1}(f)$. These properties are useful and worthwhile to study
Euler numbers and polynomials.
 } \endabstract
 \subjclass  11S80 \endsubjclass
\thanks  \endthanks
\leftheadtext{     } \rightheadtext{ }
\endtopmatter

\document

\head 1. Introduction \endhead Let $p$ be a fixed odd prime.
Throughout $\Bbb Z_p$, $\Bbb Q_p$ and $\Bbb C_p$ will respectively
denote the rings of $p$-adic integers, the fields of $p$-adic
numbers and the completion of  of the algebraic closure of $\Bbb
Q_p$. For $d$ a fixed positive integer with $(p,d)=1$, let
$$X=X_d=\varprojlim_N \Bbb Z/dp^N , \;\;X_1=\Bbb Z_p,$$
$$X^*=\bigcup\Sb 0<a<dp\\ (a,p)=1\endSb a+dp\Bbb Z_p,$$
$$a+dp^N\Bbb Z_p=\{x\in X\mid x\equiv a\pmod{dp^N}\},$$
where $a\in \Bbb Z$ lies in $0\leq a<dp^N ,$ (cf. [1], [2]).

The $p$-adic absolute value in $\Bbb C_p$ is normalized so that
$|p|_p=\frac1p .$ Let $q$ be variously considered as an
indeterminate a complex number $q \in \Bbb C$, or a $p$-adic
number $q\in\Bbb C_p .$  If $q\in\Bbb C_p,$ we assume $|q-1|_p <
1$ for $|x|_p \le 1.$ Throughout this paper, we use the following
notation :
$$[x]_q=[x:q]= \frac{1-q^x}{1-q}=1+q+q^2+\cdots+q^{x-1}.$$
We say that $f$ is uniformly differentiable function at a point
 $ a \in \Bbb Z_p$-- and denote this property by $ f \in UD(\Bbb
 Z_p )$-- if the difference quotients
 $$ F_f (x, y)=\frac{f(x)-f(y)}{x-y} ,$$ have a limit $l=
 f^{\prime}(a)$ as $( x, y) \rightarrow (a,a) ,$ cf. [1, 2].
For $f\in U D(\Bbb Z_p),$ let us start with the expression
$$\frac1{[p^N]_q} \sum_{0 \le j < p^N} q^i f(j) =\sum_{0\le j<p^N} f(j)
\mu_q(j+p^N\Bbb Z_p), \text{ cf. [1, 2], }$$ representing
$q$-analogue of Riemann sums for $f$.

The integral of $f$ on $\Bbb Z_p$ will be defined as limit
($n\rightarrow \infty$) of these sums, when it exists. An
invariant $p$-adic $q$-integral of a function $f \in {\text{UD}}
(\Bbb Z_p) $ on $\Bbb Z_p$ is defined by
$$I_{q}(f)=\int_{X}f(x)d\mu_q(x)=\int_{\Bbb Z_p} f(x) d\mu_q (x)= \lim_{N\rightarrow \infty} \frac1{[p^N]_q}
\sum_{0\le j<p^N} f(j)q^j, \text{ see [2] }.$$ Note that if $f_n
\rightarrow f$ in $UD( Z_p)$; then
$$\int_{\Bbb Z_p} f_n(x) d \mu_q(x) \rightarrow \int_{\Bbb Z_p} f(x) d
\mu_q(x).$$ The classical Euler polynomials are defined by
$$F(t, x)=\frac{2}{e^t+1}e^{xt}=\sum_{n=0}^{\infty}E_n(x)\frac{t^n}{n!}.$$
Note that $E_n=E_n(0)$ are called $n$-th Euler numbers. In [3], an
analogue of Bernoulli number is defined by
$\frac{t}{we^{t}-1}=\sum_{n=0}^{\infty}B_n(w)\frac{t^n}{n!},$
where $w$ is the element of locally constant functions' space [3].
For $f\in UD(\Bbb Z_p)$, we define a p-adic invariant integral on
$\Bbb Z_p$ as follows:
$$I_{-1}(f)=\lim_{q\rightarrow -1} I_q(f)=\int_{\Bbb
Z_p}f(x)d\mu_{-1}(x). $$
 The purpose of this paper is to give the integral equations
 related to $I_{-1}(f)$ and to investigate some properties for
 $I_{-1}(f)$. From these integral equations, we derive the
 interesting formuae related to Euler numbers.

 \head \S 2. $p$-adic invariant integrals on $\Bbb Z_p$ \endhead

For $f\in UD(\Bbb Z_p),$ the $p$-adic $q$-integral was defined by
$$I_q(f)=\lim_{N\rightarrow \infty}\frac{1}{[p^N]_q}\sum_{0\leq
x<p^N}q^xf(x)=\lim_{N\rightarrow \infty}\sum_{0\leq
x<p^N}f(x)\mu_q(x+p^N\Bbb Z_p)=\int_{\Bbb Z_p}f(x)d\mu_q(x),\tag 1
$$ representing $p$-adic $q$-analogue of Riemann integral for $f$ (see
[2]).
 In the meaning of fermionic, we now define $I_{-1}(f)$-integral as
  $$I_{-1}(f)=\lim_{q\rightarrow -1}I_q(f)=\lim_{q \rightarrow -1}\int_{\Bbb Z_p}
f(x)d\mu_{-q}(x)=\int_{\Bbb Z_p}f(x) d\mu_{-1}(x). \tag2$$ From
(2), we derive the below Eq.(3):
$$I_{-1}(f)=\int_{\Bbb Z_p}f(x)d\mu_{-1}(x)=\lim_{N\rightarrow \infty}\sum_{x=0}^{p^N-1} (-1)^xf(x) .\tag3 $$
Let $f_1(x)=f(x+1)$. Then we note that
$$I_{-1}(f_1)=-\lim_{N\rightarrow \infty}\sum_{x=0}^{p^N-1}f(x)(-1)^x
+2f(0)=-I_{-1}(f)+2f(0). $$ Therefore we obtain the below theorem:
\proclaim{ Theorem 1} For $f\in UD(\Bbb Z_p),$ we have
$$I_{-1}(f_1)+I_{-1}(f)=2f(0), \tag4$$
where $f_1(x)$ is translation with $f_1(x)=f(x+1)$.
\endproclaim
From Theorem 1, we can derive the below theorem:
 \proclaim{ Theorem 2}
For $f\in UD(\Bbb Z_p)$, $n\in\Bbb N$, we have
$$I_{-1}(f_n)+(-1)^{n-1}I_{-1}(f)=2\sum_{x=0}^{n-1}f(x)(-1)^{n-1-x}, $$
where $f_n(x)=f(x+n).$
 \endproclaim
By using Theorem 1 and Theorem 2, we can consider the new
extension of Euler numbers and polynomials. If we take
$f(x)=\lambda{^x} e^{xt},$ ($\lambda \in \Bbb Z_p$), then we have
$$\frac{2}{\lambda e^t +1}=\int_{\Bbb
Z_p}\lambda^xe^{xt}d\mu_{-1}(x). \tag5$$ Define the new extension
of Euler numbers as follows:
$$\frac{2}{\lambda e^t
+1}=\sum_{n=0}^{\infty}E_n(\lambda)\frac{t^n}{n!}. \tag6$$ By (5)
and (6), we obtain the below theorm:

\proclaim{ Theorem 3} ( Witt's formula ) For $n\in\Bbb N$, we have
$$\int_{\Bbb Z_p}\lambda^x x^n d\mu_{-1}(x)=E_n(\lambda),
\text{$\lambda \in \Bbb Z_p$}, $$ where $E_{n}(\lambda )$ are
called analog of Euler numbers.
\endproclaim
By using $I_{-1}$-integral, we define the new extension of
classical Euler polynomials as follows:
$$\int_{\Bbb Z_p}\lambda^y (x+y)^n d\mu_{-1}(y)=E_{n}(\lambda:x).
\tag7$$ From (7), we can derive the below:
$$E_n(\lambda:x)=\sum_{k=0}^n\binom nk
E_n(\lambda)x^{n-k}=f^n\sum_{a=0}^{f-1}(-1)^a\lambda^a
E_{n}(\lambda^f :\frac{x+a}{f}), \tag8$$ where $f$ is odd positive
integer.

By (8), we easily see that
$$\frac{2e^{xt}}{\lambda
e^t+1}=\int_{X}e^{(x+y)t}d\mu_{-1}(y)=\frac{2\sum_{a=0}^{f-1}(-1)^a\lambda^a
e^{(x+a)t}}{\lambda^fe^{ft}+1}. \tag9 $$
\proclaim {Remark 1} In
[3], it was known that
$$\lim_{q\rightarrow
1}\int_{X}w^ye^{(x+y)t}d\mu_q(t)=\frac{t
e^{xt}}{we^t-1}=\sum_{m=0}^{\infty}B_m(w)\frac{t^m}{m!}, \tag
9-1$$ where $B_m(w)$ are called an analogue of Bernoulli numbers.
In the viewpoint of (9-1), we gave the new extension
($E_{n}(\lambda )$) of classical Euler number.
\endproclaim

\proclaim{ Remark 2} From using multivariate $p$-adic invariant
integral for $I_{-1}(f)$, we can easily derive the Euler
polynomials of higher order as follows:

$$e^{xt}\int_{\Bbb Z_p}\cdots\int_{\Bbb Z_p}\lambda^{y_1+\cdots
+y_r}e^{(y_1+\cdots+y_r)t}d\mu_{-1}(y_1)\cdots
d\mu_{-1}(y_r)=\left( \frac{2}{\lambda e^t+1} \right)^re^{xt}.$$
From this, we can define the extension of classical Euler
polynomials of order $r$ as follows:
$$\left(\frac{2}{\lambda e^t+1}
\right)^re^{xt}=\sum_{n=0}^{\infty}E_n^{(r)}(\lambda: x).$$
\endproclaim
Let $\chi$ be the Dirichlet's character with conductor $f(=odd)\in
\Bbb N$ and let us take $f(x)=\chi(x)e^{tx}$. From Theorem 2, we
derive the below formula:
$$\int_{X}e^{tx}\chi(x)d\mu_{-1}(x)=\frac{2\sum_{a=0}^{f-1}e^{ta}(-1)^a
\chi(a)}{e^{ft}+1}.\tag10$$ The generalized Euler numbers attached
to $\chi$ were defined by
$$\frac{2\sum_{a=0}^{f-1}e^{ta}(-1)^a
\chi(a)}{e^{ft}+1}=\sum_{n=0}^{\infty}E_{n,\chi}\frac{t^n}{n!},
\text{ cf. [4] } .\tag11$$ From (10) and (11), we derive the below
Witt's formula for the generalized Euler numbers attache to $\chi$
as follows:

\proclaim{ Theorem 4} Let $f$ be an odd positive integer and let
$\chi$ be the Dirichlet's character with conductor $f$. Then we
have
$$\int_{X}x^n\chi(x) d\mu_{-1}(x)=E_{n,\chi}. $$
\endproclaim

\enddocument